\documentclass[a4paper,10pt]{article}

\usepackage[latin1]{inputenc}
\usepackage[T1]{fontenc}
\usepackage{amsmath}
\usepackage{amsfonts}
\usepackage{amstext}
\usepackage{amsthm}
\usepackage[all]{xy}
\usepackage{geometry}
\usepackage{srcltx}
\usepackage{graphics}
\usepackage{epsfig}
\usepackage{subfigure}
\usepackage{slashed}
\usepackage{color}
\usepackage{amssymb}
\usepackage{srcltx}

\makeatletter

\@addtoreset{equation}{section} \makeatother

\newtheorem{theorem}{Theorem}[section]
\newtheorem{lemma}[theorem]{Lemma}
\newtheorem{prop}[theorem]{Proposition}
\newtheorem{cor}[theorem]{Corollary}

\theoremstyle{definition}
\newtheorem{defi}[theorem]{Definition}
\newtheorem{rem}[theorem]{Remark}

\DeclareMathOperator{\im}{im}
\DeclareMathOperator{\ind}{ind}

\DeclareMathOperator{\sfl}{sf}

\DeclareMathOperator{\re}{Re}
\DeclareMathOperator{\codim}{codim}

\title{Bifurcation of critical points along\\
 gap--continuous families of subspaces}
\author{Anna Maria Candela\footnote{Supported
by Research Funds {\sl Fondi d'Ateneo 2014}.}$\ $ and Nils Waterstraat\footnote{Supported by the Berlin Mathematical School, 
the {\sl SFB 647 ``Space--Time--Matter''}, and by the program {\sl Professori Visitatori Junior 2014} of GNAMPA--INdAM.}}

\begin{document}
\date{}
\maketitle

\medskip

\begin{abstract}
We consider the restriction of twice differentiable functionals on a 
Hilbert space to families of subspaces that vary continuously with respect to the gap metric. 
We study bifurcation of branches of critical points along these families and apply our results 
to semilinear systems of ordinary differential equations.  
\end{abstract}


\section{Introduction}
Let $H$ be a real separable Hilbert space and $\mathcal{J}:H\rightarrow\mathbb{R}$ a $C^2$--functional. 
We denote the derivative of $\mathcal{J}$ at $u\in H$ by $d_u\mathcal{J}\in\mathcal{L}(H,\mathbb{R})$ 
and in what follows we assume that $d_0\mathcal{J}=0$, i.e. $0\in H$ is a critical point of $\mathcal{J}$. 
Usually, critical points of functionals $\mathcal{J}$ on Hilbert spaces $H$ are studied 
as they can be solutions of differential equations. Correspondingly, critical points of a restriction $\mathcal{J}\mid_{H'}:H'\rightarrow\mathbb{R}$ to a subspace $H'\subset H$ may yield solutions of differential equations under additional constraints.\\
In \cite{Abbondandolo} Abbondandolo and Majer studied the {\sl Grassmannian} of a Hilbert space $H$, 
i.e. the set of all closed subspaces of $H$. As there is a canonical metric on this set, which is induced by orthogonal projections, we can define paths $\{H_t\}_{t\in[a,b]}$ in it. Clearly, for each $t \in [a,b]$ the element $0\in H_t$ is a critical point of the restriction $\mathcal{J}\mid_{H_t}:H_t\rightarrow\mathbb{R}$ as $d_0\mathcal{J}=0$, and the aim of this paper is to investigate bifurcation from this branch of critical points in a sense that we will introduce below in Definition \ref{defi-bifurcation}. Our main results show the existence of bifurcation in terms of the second derivative of $\mathcal{J}$ at the critical point $0$, which are based on \cite{SFLPejsachowicz} and \cite{JacBifIch}. 
To this aim, we introduce a family of functionals $f_t:H\rightarrow\mathbb{R}$, $t\in[a,b]$, such that 
each $f_t$ involves the orthogonal projection onto the space $H_t$, and such that its critical points 
are the critical points of the restriction $\mathcal{J}\mid_{H_t}$. 
Consequently, $0\in H$ is a critical point of any $f_t:H\rightarrow\mathbb{R}$, $t\in[a,b]$, 
and by considering the second derivative $d^2_0f_t$ of $f_t$ at $0$ we can define a 
path $\{L_t\}_{t\in[a,b]}$ of bounded selfadjoint operators by the Riesz representation theorem. 
The assumptions of our theorems ensure that each $L_t$ is actually a Fredholm operator, and we prove 
that bifurcation of critical points of $f$ along $\{H_t\}_{t\in[a,b]}$ arises if the \textit{spectral flow} of 
$L:t \mapsto L_t$ does not vanish. 
Let us recall that the spectral flow is an integer valued homotopy invariant for paths of selfadjoint Fredholm operators 
that was introduced by Atiyah, Patodi and Singer in \cite{AtiyahPatodi}. Its relevance to bifurcation theory was 
discovered in \cite{SFLPejsachowicz}. 
For example, if all operators $L_t$ have a finite Morse index $\mu_{Morse}(L_t)$, then the spectral flow of $L$ is just 
the difference of the Morse indices at the endpoints, i.e. $\mu_{Morse}(L_a) - \mu_{Morse}(L_b)$. Hence a non--vanishing spectral flow of $L$ corresponds to a jump in the Morse indices of $L$, which implies bifurcation of critical points of $f$ by a well known theorem in bifurcation theory (cf. \cite[\S 8.9]{Mawhin} or also \cite[\S II.7.1]{Kielhoeffer}). However, if $\mu_{Morse}(L_t)=+\infty$ for some $t \in [a,b]$, then the spectral flow may depend on the whole path $L$ and not only on its endpoints, which makes the theory more complicated.\\
The paper is structured as follows. In Section \ref{sec2}, we introduce some preliminaries 
that we need in order to state our theorems. We recall some facts about the 
Grassmannian of a Hilbert space $H$, essentially following Abbondandolo and Majer's paper \cite{Abbondandolo}. 
However, we also state and prove a folklore result which shows that the kernels of families of surjective 
bounded operators on $H$ yield paths in the Grassmannian and which we use in the final section in our examples. 
In Section 2 we briefly recall the definition of the spectral flow from \cite{SFLPejsachowicz}. 
In the third section, we introduce the path $L$ and state our main theorems and a corollary, 
which we prove in Section \ref{sec4}. Finally, we apply our theory to a Dirichlet problem for semilinear ordinary differential operators in Section \ref{sec5}.


\section{Grassmannians and spectral flows}\label{sec2}

As before, we let $H$ be a real separable Hilbert space of infinite dimension, we denote by $\mathcal{L}(H)$ the 
Banach space of all linear bounded operators on $H$ equipped with its standard norm $\|\cdot\|$ and by $I_H\in\mathcal{L}(H)$ the identity operator. Let us recall that a {\sl Fredholm operator} $T$ on a Hilbert space $H$ is an operator $T\in\mathcal{L}(H)$ 
such that both its kernel and its cokernel are of finite dimension. We denote the open subset 
of all Fredholm operators in $\mathcal{L}(H)$ by $\Phi(H)$.


\subsection{The Grassmannian of a Hilbert space}

In this section, we recall briefly the definition and some properties of the {\sl Grassmannian} $\mathcal{G}(H)$ of $H$, 
i.e. the set of all closed linear subspaces of $H$, where we refer for more details to the comprehensive exposition 
\cite{Abbondandolo}.\\
For every $U\in\mathcal{G}(H)$, there exists a unique orthogonal projection $P_U:H\rightarrow H$ onto $U$ and the distance
\[
d(U,V):=\|P_U-P_V\|,\,\quad U,V\in\mathcal{G}(H),
\]
makes $\mathcal{G}(H)$ a complete metric space (cf. also \cite{Kato}). 
Moreover, one can show that $\mathcal{G}(H)$ is an analytic Banach manifold, and the map 
\[
V \in \mathcal{G}(H)\ \mapsto\ P_V\in\mathcal{L}(H)
\]
embeds $\mathcal{G}(H)$ analytically into $\mathcal{L}(H)$ (cf. \cite[Proposition 1.1]{Abbondandolo}). In what follows, we denote by $\{V_t\}_{t\in[a,b]}$ paths in $\mathcal{G}(H)$, i.e. continuous maps $[a,b]\rightarrow\mathcal{G}(H)$, $t\mapsto V_t$.

\begin{lemma}
The connected components of $\mathcal{G}(H)$ are the sets
\[
\mathcal{G}_{nk}(H)=\{V\in\mathcal{G}(H):\dim V=n,\,\codim V=k\},
\] 
with $n$, $k\in\mathbb{N}\cup\{+\infty\}$ such that $k+n=+\infty$.
\end{lemma}    

\begin{proof}
Let us first recall that if $\|P_U-P_V\|<1$ for $U,V\in\mathcal{G}(H)$, then $\dim U=\dim V$ and $\dim U^\perp=\dim V^\perp$ (cf. \cite[I.4.6]{Kato}). Consequently, if $U$ and $V$ belong to the same component of $\mathcal{G}(H)$, then they must have 
both the same dimension and the same codimension. \\
Now, let us assume that $U,V\in\mathcal{G}_{nk}(H)$ for some $k$, $n$ such that $k+n=+\infty$. 
Since $H$ is separable, it is easy to construct an orthogonal operator $O:H\rightarrow H$ such that $O(U)=V$. 
Denoting by $\mathcal{O}(H)$ the subspace of $\mathcal{L}(H)$ consisting of all orthogonal operators,
it is easily seen from functional calculus that $\mathcal{O}(H)$ is connected\footnote{Actually, even more is true: in \cite{Kuiper}
Kuiper proved that $\mathcal{O}(H)$ is contractible.}. 
Hence, there is a path $M:[0,1]\rightarrow\mathcal{O}(H)$ joining the identity operator $I_H$ to $O$. Finally, since $P_{M_t(U)}=M^{-1}_tP_UM_t$ for each $t\in [0,1]$, we have that $\{M_t(U)\}_{t\in [0,1]}$ is continuous and so a path in $\mathcal{G}(H)$ that joins $U$ to $V$.  
\end{proof}
\noindent
\begin{rem}
A computation of all homotopy groups $\pi_i(\mathcal{G}_{nk}(H))$, $i\in\mathbb{N}$, can be found in \cite[Section 2]{Abbondandolo}.
\end{rem}
\noindent
The following lemma is essentially well known (cf. e.g. \cite[Appendix A]{BoossZhu}), but as we are not aware of a proof in the literature, we include it here for the sake of completeness. The reader may compare it with a related assertion on Banach bundles, which can be found e.g. in \cite{BanachBundles} and also \cite{IndBundIch}, and on which our argument is based. 

\begin{lemma}\label{surjops}
Let $A:[a,b]\rightarrow\mathcal{L}(H,X)$ be a continuous family of bounded surjective operators, 
where $X$ is a Banach space and $\mathcal{L}(H,X)$ denotes the Banach space of all bounded linear operators. Then
\[
\{\ker A_t\}_{t\in[a,b]}:=\{u\in H:\ A_t u=0\}_{t\in [a,b]}
\]
is a path in $\mathcal{G}_{nk}(H)$, where $k=\dim X$ and $n=\dim H-\dim X$.
\end{lemma}

\begin{proof}
Let us first fix some $t_0\in[a,b]$. Since $A_{t_0}$ is surjective, there exists $M_0\in\mathcal{L}(X,H)$ 
such that $A_{t_0}M_0=I_X$, with $I_X$ the identity operator on $X$. 
From the fact that the invertible elements in $\mathcal{L}(X)$ are open, we see 
that $A_t M_0$ is invertible for all $t$ in a neighbourhood $\mathcal{I}_0$ of $t_0$.\\
Now, if we set $M_{0,t}:=M_0(A_t M_0)^{-1}$ for $t\in \mathcal{I}_0$, then $A_t M_{0,t}=I_X$.\\ 
Note that if $M_1,M_2\in\mathcal{L}(X,H)$ are such that $A_tM_i=I_X$, then $A_t(\alpha M_1+(1-\alpha)M_2)=I_X$ 
for all $0\leq\alpha\leq 1$. Consequently, by using a partition of unity, we may conclude that 
there exists a path $M:[a,b]\rightarrow\mathcal{L}(X,H)$ such that $A_t M_t=I_X$ for all $t\in [a,b]$.\\
Defining $R_t:=M_t A_t\in\mathcal{L}(H)$, we note that $R_t$ is a projection since 
\[
R^2_t\ =\ M_t A_t M_t A_t\ =\ M_t A_t\ =\ R_t. 
\]
Moreover, since $M_t$ is clearly injective, we infer that 
\[
\ker(R_t)=\ker(M_t A_t)=\ker(A_t)
\]
so that $Q_t:=I_H-R_t$ is a continuous family of projections such that $\im(Q_t)=\ker(A_t)$. 
Thus, taking 
\[
P_t\ =\ Q_t Q^\ast_t(Q_t Q^\ast_t+(I_{H}-Q^\ast_t)(I_{H}-Q_t))^{-1},
\]
it follows by \cite[Lemma 12.8 a)]{BoossBuch} that
$\{P_t\}_{t\in[a,b]}$ is a continuous family of orthogonal projections such that $\im(P_t)=\ker(A_t)$. Hence, $\{\ker(A_t)\}_{t\in [a,b]}$ is a continuous family of subspaces in $\mathcal{G}(H)$.\\
Finally, that $\ker(A_t)\in\mathcal{G}_{nk}(H)$ with $k=\dim X$ and $n=\dim H-\dim X$ is an immediate consequence of the 
rank--nullity theorem in linear algebra.    
\end{proof}


\subsection{The spectral flow}\label{section-bif}
We denote by $\Phi_S(H)\subset\Phi(H)$ the subspace of all selfadjoint Fredholm operators, 
which is well known to consist of three connected components (cf. \cite{AtiyahSinger}). 
Two of them are given by
\[ \begin{split}
&\Phi^+_S(H)=\{L\in \Phi_S(H):\sigma_{ess}(L)\subset(0,+\infty)\},\\
&\Phi^-_S(H)=\{L\in \Phi_S(H):\sigma_{ess}(L)\subset(-\infty,0)\},
\end{split}
\]
where $\sigma_{ess}(L)=\{\lambda\in\mathbb{R}:\, L-\lambda I_H\notin\Phi_S(H)\}$ is the {\sl essential spectrum} of an operator $L\in\Phi_S(H)$. 
Their elements are called {\sl essentially positive} or {\sl essentially negative},
respectively, and it is readily seen that both of these spaces are contractible. 
Elements of the remaining component $\Phi^i_S(H)=\Phi_S(H)\setminus(\Phi^+_S(H)\cup\Phi^-_S(H))$ are called {\sl strongly indefinite}, 
and in contrast to $\Phi^+_S(H)$ and $\Phi^-_S(H)$, this space has a non--trivial topology. 
Indeed, $\Phi^i_S(H)$ has the same homotopy groups as the stable orthogonal group (cf. \cite{Schroeder}) and the spectral flow
provides an explicit isomorphism between its fundamental group and the integers. 
There are several different, but equivalent, constructions of the spectral flow in the literature. 
Here, we follow the approach developed by
Fitzpatrick, Pejsachowicz and Recht in \cite{SFLPejsachowicz}, and we refer to the introduction of 
\cite{JacBifIch} for further references on the subject.\\
We call two selfadjoint invertible operators in $\mathcal{L}(H)$ {\sl Calkin equivalent} if $S-T$ is compact. 
It is well known that in this case the {\sl relative Morse index}
\[
\mu_{rel}(S,T)=\dim(E^-(S)\cap E^+(T))-\dim(E^+(S)\cap E^-(T))
\]
is well defined and finite, where $E^-(\cdot)$ and $E^+(\cdot)$ denote the negative and positive subspaces 
of a selfadjoint operator for which $0$ is an isolated point of the spectrum.\\
From the second resolvent identity it follows that for Calkin equivalent operators $S,T$, also the difference of the associated resolvent operators
\[
(\lambda-T)^{-1}-(\lambda-S)^{-1}=(\lambda-T)^{-1}(T-S)(\lambda-S)^{-1},\qquad\lambda\notin\sigma(T)\cup\sigma(S),
\]
is compact whenever it is defined, where $\sigma(T)$ and $\sigma(S)$ denote the spectrum of $T$ and $S$, respectively. 
Finally, since the set of compact operators is closed in $\mathcal{L}(H)$, it follows that also the difference of the spectral projections
\[
P_{[a,b]}(T)-P_{[a,b]}(S)\ =\ \re\left(\frac{1}{2\pi i}\int_{\Gamma}\big[(\lambda-T^\mathbb{C})^{-1}-
(\lambda-S^\mathbb{C})^{-1}\big]\ d\lambda\right)
\]
is compact, where $a,b$ do not belong to $\sigma(S)\cup\sigma(T)$ and $\Gamma$ is the circle around $\frac{a+b}{2}$ in $\mathbb{C}$ intersecting the real axis at $a$ and $b$. Here, $S^\mathbb{C}$ and $T^\mathbb{C}$ denote the complexification 
of operators and $\re$ the real part of an operator on a complexified Hilbert space (cf. \cite[Subsection 2.1]{domainshrinking}
for more details).\\
The group $GL(H)$ of all invertible operators on $H$ acts on $\Phi_S(H)$ by mapping $M\in GL(H)$ and $L\in\Phi_S(H)$ to $M^\ast LM$, which is called the \textit{cogredient action}. One of the main theorems in \cite{SFLPejsachowicz} 
states that for any path $L:[a,b]\rightarrow\Phi_S(H)$ there exist a path $M:[a,b]\rightarrow GL(H)$ and a selfadjoint invertible 
operator $J\in\Phi_S(H)$, such that $M^\ast_tL_tM_t=J+K_t$ with $K_t$ selfadjoint and compact for each $t\in[a,b]$.

\begin{defi}
Let $L:[a,b]\rightarrow\Phi_S(H)$ be a path such that $L_a$ and $L_b$ are invertible. 
The \textit{spectral flow} of $L$ is the integer
\[
\sfl(L,[a,b])\ =\ \mu_{rel}(J+K_a,J+K_b),
\]
where $J+K:[a,b]\rightarrow\Phi_S(H)$ is any path of compact selfadjoint perturbations $K_t$, $t\in[a,b]$, of a selfadjoint invertible operator $J\in\Phi_S(H)$ which is cogredient to $L$.
\end{defi}
\noindent
It follows from well known properties of the relative Morse index that the spectral flow does not 
depend on the choice of the path $J+K$, and moreover it has the following properties:
\begin{itemize}
\item[{\sl (i)}] If $L_{t}$ is invertible for all $t\in [a,b]$, then $\sfl(L,[a,b])=0$.
\item[{\sl (ii)}] If $H_1$ and $H_2$ are separable Hilbert spaces and the paths $L_1:[a,b]\rightarrow \Phi_S(H_1)$ 
	and $L_2:[a,b]\rightarrow \Phi_S(H_2)$ have invertible endpoints, then
\[
\sfl(L_1\oplus L_2,[a,b])\ =\ \sfl(L_1,[a,b])+\sfl(L_2,[a,b]).
\]
\item[{\sl (iii)}] Let $h:[0,1]\times [a,b]\rightarrow\Phi_S(H)$ be a homotopy such that
$h(s,a)$ and $h(s,b)$ are invertible for all $s\in [0,1]$. Then,
	\[
	\sfl(h(0,\cdot),[a,b])\ =\ \sfl(h(1,\cdot),[a,b]).
	\]
  \item[{\sl (iv)}] If $L_t\in\Phi_S^+(H)$, $t\in [a,b]$, and $L_a$, $L_b$ are invertible, 
  then the spectral flow of $L$ is the difference of the Morse indices at its endpoints:
	\[
	\sfl(L,[a,b])\ =\ \mu_{Morse}(L_a)-\mu_{Morse}(L_b),
	\]
	where
	\begin{equation}\label{star}
	\mu_{Morse}(L_t)\ =\ \sup\dim\{V\subset H:\, \langle L_t u,u\rangle_H<0\,\,\, \hbox{for all $u\in V\setminus\{0\}$}\}.
	\end{equation}
\end{itemize}
\noindent
Finally, let us note that the spectral flow is actually uniquely characterised 
by the properties {\sl (i)}--{\sl (iv)} above (cf. \cite{UnSfl}). A further uniqueness theorem for the spectral flow, 
which is based on the different but equivalent construction \cite{Phillips}, can be found in \cite[\S 5.2]{LeschSpecFlowUniqu}.


\section{Bifurcation along gap continuous paths of subspaces}\label{sec3}

As before, let $H$ be a real Hilbert space and $\mathcal{J}:H\rightarrow\mathbb{R}$ a 
$C^2$-functional having $0$ as a critical point. We denote by $d_u\mathcal{J}\in\mathcal{L}(H,\mathbb{R})$ 
the derivative of $\mathcal{J}$ at $u\in H$, and we let $T$ be the Riez representation of the Hessian 
$d^2_0\mathcal{J}:H\times H\rightarrow\mathbb{R}$ of $\mathcal{J}$ at $0$, i.e. 
the unique selfadjoint operator $T\in\mathcal{L}(H)$ which satisfies
\begin{align}\label{T}
d^2_0\mathcal{J}[u,v]\ =\ \langle Tu,v\rangle_H,\quad u,v\in H.
\end{align}
Let $\{H_t\}_{t\in [a,b]}\subset\mathcal{G}(H)$ be a gap continuous path of closed subspaces of $H$ 
for some real numbers $a<b$, and let us point out that $0\in H$ is in any $H_t$, $t\in [a,b]$. 
In what follows we denote by $\mathcal{J}\mid_{H_t}:H_t\rightarrow\mathbb{R}$ the restriction of the 
functional $\mathcal{J}$ to the closed subspace $H_t\subset H$. Note that $0\in H$ is a critical point 
of all $\mathcal{J}\mid_{H_t}$, $t\in[a,b]$, which is a direct consequence of the uniqueness of the derivative.

\begin{defi}\label{defi-bifurcation}
We say that $t^\ast\in[a,b]$ is a {\sl bifurcation point of $\mathcal{J}$ along $\{H_t\}_{t\in[a,b]}$} 
if there exist two sequences $(t_n)_{n}\subset [a,b]$ and 
$(u_n)_{n}\subset H$ such that 
\begin{enumerate}
	\item[{\sl (i)}] $t_n\rightarrow t^\ast$ in $[a,b]$ and
	 $u_n\rightarrow 0$ in $H$ as $n\rightarrow+\infty$;
		\item[{\sl (ii)}] $u_n\in H_{t_n}$ and $u_n\neq 0$ for all $n\in\mathbb{N}$;
	\item[{\sl (iii)}] $u_n$ is a critical point of $\mathcal{J}\mid_{H_{t_n}}$ for all $n\in\mathbb{N}$. 
\end{enumerate}
\end{defi}  
\noindent
Since $\{H_t\}_{t\in [a,b]}$ is a continuous path of subspaces, there exists a family $P_t$, 
$t\in [a,b]$, of orthogonal projections such that $\im P_t=H_t$. We set $P^\perp_t:=I_H-P_t$, and define
\begin{align}\label{L}
L_t=P_tTP_t+P^\perp_t \quad \hbox{for each $t\in [a,b]$,}
\end{align}
which is a continuous path of selfadjoint operators in $\mathcal{L}(H)$. We call $\{H_t\}_{t\in [a,b]}$ \textit{admissible} if
both operators
\[
P_a TP_a:H_a\rightarrow H_a\quad \hbox{and}\quad P_b TP_b:H_b\rightarrow H_b
\]
are invertible. Since $H_t$ and $H^\perp_t$ reduce $L_t$, and $L_t\mid_{H^\perp_t}=I_{H^\perp_t}$ is invertible, 
we see at once that $L_a$ and $L_b$ are invertible if $\{H_t\}_{t\in [a,b]}$ is admissible.\\
Now, let us state our main theorems and a corollary, which we are proving in the next section. 

\begin{theorem}\label{thm}
Let $\{H_t\}_{t\in [a,b]}$ be an admissible path in $\mathcal{G}_{nk}(H)$ such that either 
$n\neq +\infty$ or $k\neq +\infty$, and let us assume that the operator $T$ introduced in \eqref{T} is Fredholm.\\
Then the operators $L_t$ in \eqref{L} are Fredholm, and if $\sfl(L,[a,b])\neq 0$, then there is a bifurcation point of $\mathcal{J}$ along $\{H_t\}_{t\in [a,b]}$. Moreover, if $n\neq +\infty$ and 
$\{H_t\}_{t\in[a,b]}$ is analytic, then there are at least
\begin{equation}\label{integer}
\left\lfloor\frac{|\sfl(L,[a,b])|}{n}\right\rfloor
\end{equation}
distinct bifurcation points (here, $\lfloor\cdot\rfloor$ denotes the integer part of a positive real number).
\end{theorem}
\noindent
Note that the case in which the path $\{H_t\}_{t\in [a,b]}$ is in the connected component $\mathcal{G}_{\infty,\infty}(H)$ of $\mathcal{G}(H)$ is excluded in Theorem \ref{thm}.  
Our second theorem deals with this setting, but we have to impose a restriction on the form of the operator $T$.

\begin{theorem}\label{thmII}
We assume that $T=I_H+K$ for some compact operator $K$, and that $\{H_t\}_{t\in [a,b]}$ is an admissible path in $\mathcal{G}_{\infty,\infty}(H)$. Then the operators $L_t$ in \eqref{L} are Fredholm, and if $\sfl(L,[a,b])\neq 0$, then there is a bifurcation point of $\mathcal{J}$ along $\{H_t\}_{t\in [a,b]}$.
\end{theorem}
\noindent
Let us point out that $L_t\in\Phi^+_S(H)$, $t\in[a,b]$, and so 
\[
\begin{split}
\sfl(L,[a,b])&=\mu_{Morse}(L_a)-\mu_{Morse}(L_b)=\mu_{Morse}(T\mid_{H_a})-\mu_{Morse}(T\mid_{H_b}),
\end{split}
\]
in each of the following cases:
\begin{itemize}
\item if $n\neq +\infty$ in Theorem \ref{thm}, since each $L_t$ is positive on the subspace $H^\perp_t$ which is of finite codimension;
\item if $T\in\Phi^+_S(H)$ in Theorem \ref{thm}, as $\mu_{Morse}(L_t)\leq\mu_{Morse}(T)$ for all $t\in [a,b]$;
\item for all compact operators $K$ in Theorem \ref{thmII} by the same argument as in the previous item. 
\end{itemize}
Finally, we will prove in the subsequent section a corollary of the proof of Theorem \ref{thm}, which rephrases a well known 
fact from bifurcation theory in our setting. Let us point out that both Theorem \ref{thm} and 
Theorem \ref{thmII} do not give any information about the location of the bifurcation point in the interval $(a,b)$. 

\begin{cor}\label{cor}
We assume that either the assumptions of Theorem \ref{thm} or the ones of Theorem \ref{thmII} hold. 
If $t^\ast$ is a bifurcation point, then 
\[
\im(T\mid_{H_{t^\ast}})\cap H^\perp_{t^\ast}\neq\{0\}.
\]
\end{cor}
\noindent


\section{Proofs of the main theorems}\label{sec4}

Our proofs are based on the main theorem of \cite{JacBifIch}, which deals with the relation between the spectral flow and the bifurcation theory that was previously established in \cite{SFLPejsachowicz}. Let us first briefly recall this theorem: We assume that $f:[a,b]\times H\rightarrow\mathbb{R}$ is a continuous map such that each $f_t:=f(t,\cdot)$ is $C^2$ and all its derivatives depend continuously on $t\in[a,b]$. In what follows, if $0\in H$ is a critical point of all $f_t$, we call $t^\ast$ a {\sl bifurcation point of critical points of the functional $f$} if there exist two sequences 
$(t_n)_n \subset [a,b]$ and $(u_n)_n \subset H \setminus\{0\}$ such that 
$t_n\rightarrow t^\ast$ in $[a,b]$, $u_n\rightarrow 0$ in $H$ and $u_n$ is a critical point of $f_{t_n}$
for all $n\in\mathbb{N}$.
The second derivatives $d^2_0f_t$ of $f_t$, $t\in [a,b]$, define selfadjoint operators $L_t$ by the Riesz representation theorem, i.e.  

\begin{align*}
d^2_0f_t[u,v]\ =\langle L_tu,v\rangle_H,\quad u,v\in H,\quad t\in[a,b].
\end{align*}
The following theorem is the main result of \cite{JacBifIch} (cf. also \cite{Fitzpatrick}): 

\begin{theorem}\label{theorem}
If each $L_t$, $t\in [a,b]$, is a Fredholm operator, both $L_a$ and $L_b$ are invertible and $\sfl(L,[a,b])\neq 0$, 
then there is a bifurcation point of critical points of the functional $f$ in $(a,b)$. Moreover, if there are only 
finitely many $t\in (a,b)$ such that $\ker(L_t)\neq 0$ and
\[
m\ :=\ \sup_{t\in (a,b)}\dim\ker(L_t)\ <\ +\infty,
\]
then the number of bifurcation points is at least
\[
\left\lfloor\frac{|\sfl(L,[a,b])|}{m}\right\rfloor.
\]
\end{theorem}
\noindent  
Now, in the setting of Section \ref{sec3}, we define a one--parameter family of functionals by
\[
f_t : u \in H\ \mapsto\ f_t(u)=\mathcal{J}(P_t u)+\frac{1}{2}\|P^\perp_t u\|^2 \in \mathbb{R}.
\]

\begin{lemma}\label{equivalence}
The critical points of $f_t$ are precisely the critical points of $\mathcal{J}\mid_{H_t}$, $t\in[a,b]$.
\end{lemma}

\begin{proof}
If $u$ is a critical point of $f_t$, then
\begin{align}\label{derivative}
d_uf_t(v)\ =\ d_{P_t u}\mathcal{J}(P_t v)+\langle P^\perp_t u,P^\perp_t v\rangle\ =\ 0\qquad 
\hbox{for all $v\in H$.}
\end{align}
In particular, taking $v=P^\perp_tu$, it follows that
\[
0\ =\ d_{P_t u}\mathcal{J}(P_t P^\perp_t u)+\|P^\perp_t u\|^2\ = \ d_{P_t u}\mathcal{J}(0)+\|P^\perp_t u\|^2
\]
as $P_t P^\perp_t u=0$. Hence, $P^\perp_t u=0$ and we see that $u\in H_t$. Consequently, we obtain from \eqref{derivative} that 
\[
0\ =\ d_{P_t u}\mathcal{J}(P_t v)\ =\ d_u\mathcal{J}(v) \qquad \hbox{for all $v\in H_t$,}
\]
which shows that $u$ is a critical point of the restriction of $\mathcal{J}$ to $H_t$.\\
Conversely, if $u$ is a critical point of the restriction of $\mathcal{J}$ to $H_t$, then $u\in H_t$ and
\[
d_uf_t(v)\ =\ d_{P_t u}\mathcal{J}(P_t v)+\langle P^\perp_t u,P^\perp_t v\rangle\ =\ d_u\mathcal{J}(P_t v)
\]
which vanishes for all $v\in H$ as $P_t v\in H_t$.
\end{proof}
\noindent
Consequently, it follows from Definition \ref{defi-bifurcation} and Lemma \ref{equivalence} that $t^\ast\in [a,b]$ 
is a bifurcation point of $\mathcal{J}$ along $\{H_t\}_{t\in [a,b]}$ if and only if 
it is a bifurcation point for the family of functionals $f_t$.\\ 
By applying Theorem \ref{theorem}, for each $t \in [a,b]$ we have to consider the Hessian of $f_t$ at the critical point $0\in H$, 
which is given by
\[
d^2_0f_t[u,v]\ =\ d^2_0\mathcal{J}[P_t u,P_t v]+\langle P^\perp_t u,P^\perp_t v\rangle\qquad 
\hbox{for all $u,v\in H$.}
\]
Using that $P^\ast_t=P_t$ and $(P^\perp_t)^\ast=(P^\perp_t)^2=P^\perp_t$, we see that the corresponding Riesz representation is given by 
\[
L_t=P_t TP_t+P^\perp_t.
\]
Note that these are exactly the operators introduced in \eqref{L}.\\
Now, we deduce Theorems \ref{thm} and \ref{thmII} from Theorem \ref{theorem} but before we note for later reference the following immediate consequence of the definition of Fredholm operators.

\begin{lemma}\label{Fredholm}
If $H_1$, $H_2$ are Hilbert spaces and $T_1:H_1\rightarrow H_1$, $T_2:H_2\rightarrow H_2$ are Fredholm operators, then
\[
T_1\oplus T_2: u_1+u_2 \in H_1\oplus H_2\ \mapsto\ (T_1\oplus T_2)(u_1+u_2) = T_1u_1+T_2u_2 \in H_1\oplus H_2
\]
is a Fredholm operator of index $\ind(T_1\oplus T_2)=\ind(T_1)+\ind(T_2)$.
\end{lemma}
\noindent
In what follows, we will apply Lemma \ref{Fredholm} to $L_t\mid_{H_t}:H_t\rightarrow H_t$ and $L_t\mid_{H^\perp_t}:H^\perp_t\rightarrow H^\perp_t$.

\begin{proof}[Proof of Theorem \ref{thm}]
Let us first assume that $n\neq +\infty$. Then, by Lemma \ref{Fredholm} the operator $L_t$ is Fredholm as it is invertible on the subspace $H^\perp_t$ and Fredholm on the finite dimensional space $H_t$. 
Furthermore, $L_a$ and $L_b$ are invertible by assumption and so Theorem \ref{thm} follows from Theorem \ref{theorem}. This shows the first part of the assertion of Theorem \ref{thm}. Now, if $\{H_t\}_{t\in[a,b]}$ is analytic, then $P_t$ and so $L_t$ depends analytically on $t$. As in \cite[Section 2]{JacBifIch}, this implies that the set of all $t$ such that $\ker(L_t)\neq\{0\}$ is discrete. Moreover, it is readily seen that 
\begin{align}\label{kernel}
\ker L_t\ =\ \im(T\mid_{H_t})\cap H^\perp_t
\end{align}
for any $t\in[a,b]$, and so 
\[
\dim \ker(L_t)\ \le\ \dim\im(T\mid_{H_t})\ \leq\ \dim H_t\ =\ n.
\]
Hence, also \eqref{integer} follows from Theorem \ref{theorem}.\\
Let us now assume that $k\neq +\infty$. Since $L_a$ and $L_b$ are again invertible by assumption, in order to apply Theorem \ref{theorem} it is enough to show that $L_t$ is Fredholm for all $t\in(a,b)$. However, by Lemma \ref{Fredholm} we just need to prove that $P_tTP_t$ is Fredholm on $H_t$. Now the kernel and cokernel of the projection $P_t$ are $H^\perp_t$, which is of finite dimension $k< +\infty$, and so $P_t$ is a Fredholm operator. This shows that indeed $P_t TP_t$ is Fredholm as the composition of Fredholm operators is again Fredholm (cf. \cite[Theorem 3.2]{GohbergClasses}).
\end{proof}
\noindent
\begin{proof}[Proof of Theorem \ref{thmII}]
Our aim is again to apply Theorem \ref{theorem}, for which we need to prove 
that $L_t$ is Fredholm for all $t\in [a,b]$. 
However, as $k = n = +\infty$, none of the arguments used in the proof of Theorem \ref{thm} 
can be applied here. Instead, by the assumption that $T$ is a compact perturbation of the identity,
we see that
\[
L_t\ =\ P_t TP_t+P^\perp_t\ =\ P_t (I_H+K)P_t+P^\perp_t\ =\ P_t+P_t KP_t+P^\perp_t\ =\ I_H+P_t KP_t,
\] 
which is a compact perturbation of $I_H$ as the set of compact operators is an ideal in $\mathcal{L}(H)$. 
Now, $L_t$ is Fredholm by a classical result of Riesz and Schauder saying that compact perturbations of the identity are Fredholm operators (cf. \cite[Corollary XII.2.5]{GohbergClasses}).
\end{proof}


\noindent
\begin{proof}[Proof of Corollary \ref{cor}]
We have already shown that a bifurcation point $t^* \in (a,b)$ exists under the assumptions of Theorem \ref{thm} or Theorem \ref{thmII}, respectively. We now argue by contradiction and assume that $\im(T\mid_{H_{t^\ast}})\cap H^\perp_{t^\ast}= \{0\}$. Then $\ker(L_{t^\ast})=\{0\}$ by \eqref{kernel} and so $L_{t^\ast}$ is invertible as it is Fredholm of index $0$.\\
We now consider the map
\[
F:(t,u) \in [a,b]\times H\ \mapsto\ F(t,u)=d_uf_t \in \mathcal{L}(H,\mathbb{R})
\] 
and we note that $F(t,0)=0$ for all $t\in [a,b]$ by assumption. Since $d_0F_{t^\ast}(u)[v]=\langle L_{t^\ast}u,v\rangle$, $u,v\in H$, and as $L_{t^\ast}$ is invertible, we see that 
$d_0F_{t^\ast}:H\rightarrow\mathcal{L}(H,\mathbb{R})$ is invertible. Consequently, by the implicit function theorem all solutions of the equation $F(t,u)=0$ in a neighbourhood of $(t^\ast,0)\in [a,b]\times H$ 
are of the form $(t,0)$ and so $t^\ast$ is not a bifurcation point of critical points of $f_t$. This is a contradiction, as the bifurcation points of critical points of $f_t$ are the bifurcation points of $\mathcal{J}$ along $\{H_t\}_{t\in [a,b]}$. 
\end{proof}


\section{An example}\label{sec5}
Throughout this section, we set $I:=[0,1]$ and we
denote by $H^1_0(I,\mathbb{R}^n)$ the Hilbert space of all absolutely continuous functions 
$u:I\rightarrow\mathbb{R}^n$ such that the derivative $u'$ is square integrable.\\
Our aim is to investigate the existence of nontrivial solutions 
for the semilinear system of ordinary differential equations 
\begin{equation}\label{ODE}
\left\{
\begin{array}{ll}
-(A(x)u'(x))'+\nabla_\xi g(x,u(x))\ =\ 0, &x\in I,\\
u(0)=u(1)=0, &
\end{array}
\right.
\end{equation} 
where $A:I\rightarrow GLS(n,\mathbb{R})$ is a smooth family of invertible symmetric matrices, 
and $g:I\times\mathbb{R}^n\rightarrow\mathbb{R}$, $g=g(x,\xi)$, is a $C^2$ function 
such that $\nabla_\xi g(x,0)=0$ for all $x\in I$.\\
Let us consider the functional
$\mathcal{J}:H^1_0(I,\mathbb{R}^n)\rightarrow\mathbb{R}$ such that
\[
\mathcal{J}(u)\ =\ \frac{1}{2}\int^1_0 \langle A(x)u'(x),u'(x)\rangle\, dx\ +\ \int^1_0 g(x,u(x))\,dx.
\]
It is well known (see, e.g., \cite[Proposition B.34]{Rabinowitz}) that $\mathcal{J}$ is of class $C^2$ 
in $H^1_0(I,\mathbb{R}^n)$ and
\begin{equation}\label{eqdiff}
d_u\mathcal{J}(v)\ =\ \int^1_0{\langle A(x)u'(x),v'(x)\rangle\, dx}+\int^1_0{\langle\nabla_\xi g(x,u(x)),v(x)\rangle\, dx}
\end{equation}
for any $u, v\in H^1_0(I,\mathbb{R}^n)$. Hence the critical points of $\mathcal{J}$ are precisely the 
weak solutions of problem \eqref{ODE}.\\ 
In particular, $0\in H^1_0(I,\mathbb{R}^n)$ is a critical point and one can show that the corresponding Hessian is given by
\[
d^2_0\mathcal{J}[u,v]\ =\ \int^1_0{\langle A(x)u'(x),v'(x)\rangle\, dx}+\int^1_0{\langle S(x)u(x),v(x)\rangle\, dx}
\quad \hbox{for all $u,v\in H^1_0(I,\mathbb{R}^n)$,}
\]
where $S(x)=D^2_\xi g(x,0)$ is a family of symmetric matrices which is continuous with respect to $x$.\\
Let us recall that for every $t\in I$ there is the evaluation map 
\[
ev_t: u \in H^1_0(I,\mathbb{R}^n)\ \mapsto\ ev_t(u)=u(t) \in \mathbb{R}^n,
\]
which is a bounded linear operator that is surjective if $t\in(0,1)$. Moreover, $ev_t$ depends continuously on $t$ in $(0,1)$. 
Indeed, for every $t_0\in (0,1)$ and $u\in H^1_0(I,\mathbb{R}^n)$, we have
\[
u(t)\ =\ u(t_0)+\int^t_{t_0}{u'(s)\,ds},\quad t\in I,
\]
which implies that
\[
\|ev_t-ev_{t_0}\|\ \leq\ \sqrt{|t-t_0|}.
\]
Now, Lemma \ref{surjops} shows that we get for every $0<a<b<1$ a continuous family of 
subspaces $\{H_t\}_{t\in[a,b]}$ by
\[
H_t\ =\ \ker(ev_t)\ =\ \{u\in H^1_0(I,\mathbb{R}^n):\, u(t)=0\},
\]
and moreover, it follows by a straightforward computation that the orthogonal projection in $H^1_0(I,\mathbb{R}^n)$ onto $H_t$ is given by

\begin{align}\label{ortproj}
(P_tu)(x)=u(x)-\frac{\min\{t,x\}-tx}{t(1-t)}\,u(t),\quad x\in I.
\end{align}
\noindent
\begin{defi}\label{defpoint}
We say that $t^\ast\in(a,b)$ is a {\sl bifurcation point for \eqref{ODE}}
if there exist two sequences $(t_k)_k \subset [a,b]$ and $(u_k)_k \subset H^1_0(I,\mathbb{R}^n)$
such that
\begin{itemize}
\item[{\sl (i)}]
$t_k\rightarrow t^\ast$ in $[a,b]$ and $u_k\rightarrow 0$ in $H^1_0(I,\mathbb{R}^n)$ as $k \to +\infty$; 
\item[{\sl (ii)}] $u_k\not\equiv 0$ for each $k\in \mathbb{N}$; 
\item[{\sl (iii)}] for every $k\in \mathbb{N}$, the restriction 
$u_{0,k}:=u_k\mid_{[0,t_k]}$ satisfies 
		\[
		-(A(x)u'_{0,k}(x))'+\nabla_\xi g(x,u_{0,k}(x))=0,\qquad x\in[0,t_k];
		\]
\item[{\sl (iv)}] for every $k\in \mathbb{N}$, the restriction $u_{1,k}:=u_k\mid_{[t_k,1]}$ satisfies 
		\[
		-(A(x)u'_{1,k}(x))'+\nabla_\xi g(x,u_{1,k}(x))=0, \qquad x\in[t_k,1],
		\]
	\item[{\sl (v)}] $u_{0,k}(t_k)=u_{1,k}(t_k)=0$ for each $k\in\mathbb{N}$. 
\end{itemize}
\end{defi}
\noindent
Let us note that the two restrictions $u_{0,k}$ and $u_{1,k}$ in Definition \ref{defpoint} define a global solution of \eqref{ODE} if and only if $u'_{0,k}(t_k)=u'_{1,k}(t_k)$. 

\begin{lemma}\label{bif}
There is a bifurcation point of \eqref{ODE} at $t^\ast\in(a,b)$ if and only if $t^\ast$ is a bifurcation point of $\mathcal{J}$ along $\{H_t\}_{t\in[a,b]}$.
\end{lemma}

\begin{proof}
If $t^\ast\in(a,b)$ is a bifurcation point of \eqref{ODE}, then there are sequences
$(t_k)_k \subset [a,b]$ and $(u_k)_k \subset H^1_0(I,\mathbb{R}^n)$ which satisfy the properties 
{\sl (i)}--{\sl (v)} in Definition \ref{defpoint}. Hence, for all $v\in H_{t_k}$ we have that
\[
\int^{t_k}_0{\langle A(x)u'_{0,k}(x),v'(x)\rangle\, dx}+\int^{t_k}_0{\langle\nabla_\xi g(x,u_{0,k}(x)),v(x)\rangle\, dx}\ =\ 0
\]
and
\[
\int^1_{t_k}{\langle A(x)u'_{1,k}(x),v'(x)\rangle\, dx}+\int^1_{t_k}{\langle\nabla_\xi g(x,u_{1,k}(x)),v(x)\rangle\, dx}\ =\ 0.
\]
It follows by \eqref{eqdiff} that $u_k\in H^1_0(I,\mathbb{R}^n)$ is a non--trivial critical point 
of $\mathcal{J}\mid_{H_{t_k}}$, and as $u_k\rightarrow 0$, we see that $t^\ast$ is a bifurcation point of $\mathcal{J}$ along $\{H_t\}_{t\in[a,b]}$ (see Definition \ref{defi-bifurcation}).\\
Conversely, let $(t_k)_k \subset [a,b]$ and $(u_k)_k \subset H^1_0(I,\mathbb{R}^n)\setminus \{0\}$
be such that $u_k\in H_{t_k}$ is a critical point of $\mathcal{J}\mid_{H_{t_k}}$, 
with $t_k\rightarrow t^\ast$ and $u_k\rightarrow 0$ in $H^1_0(I,\mathbb{R}^n)$. 
Setting $u_{0,k}$ and $u_{1,k}$ as in {\sl (iii)} and {\sl (iv)} of Definition \ref{defpoint}, we get that   
\[
\int^{t_k}_0{\langle A(x)u'_{0,k}(x),v'(x)\rangle\, dx}+\int^{t_k}_0{\langle\nabla_\xi g(x,u_{0,k}(x)),v(x)\rangle\, dx}\ =\ 0\quad 
\hbox{for all $v\in H^1_0([0,t_k],\mathbb{R}^n)$}
\]
and
\[
\int^1_{t_k}{\langle A(x)u'_{1,k}(x),v'(x)\rangle\, dx}+\int^1_{t_k}{\langle\nabla_\xi g(x,u_{1,k}(x)),v(x)\rangle\, dx}\ =\ 0\quad 
\hbox{for all $v\in H^1_0([t_k,1],\mathbb{R}^n)$.}
\]
Hence $u_k$ satisfies {\sl (iii)} and {\sl (iv)} in Definition \ref{defpoint},
while {\sl (v)} is an immediate consequence of the definition of $H_{t_k}$. Thus $t^\ast$ is a bifurcation point of \eqref{ODE}.
\end{proof}
 \noindent
As by Lemma \ref{bif} the existence of bifurcation points of \eqref{ODE} can be reduced to the study of bifurcation
points of the functional $\mathcal{J}$ on $\{H_t\}_{t\in[a,b]}$, we now want to assume that 
the bilinear form $d^2_0\mathcal{J}$ is non--degenerate on $H_a$ and $H_b$.
This implies that the path $\{H_t\}_{t\in[a,b]}$ is admissible. Moreover, a straightforward computation shows that  that the operator $T:H^1_0(I,\mathbb{R}^n)\rightarrow H^1_0(I,\mathbb{R}^n)$ from \eqref{T} is given in this case by 
\[\begin{split}
Tu(x)\ =\ &\int_0^x A(s)u'(s)\ ds - x \int_0^1 A(s)u'(s)\ ds \\
&- \int^x_0{\int^s_0{S(\tau)u(\tau)\ d\tau}ds} + x \int^1_0{\int^s_0{S(\tau)u(\tau)\ d\tau}ds}.
\end{split}
\]
Consequently, by using \eqref{ortproj} we can write down the path $L_t$ in \eqref{L} explicitely and so we have everything at hand in order to claim the existence of a bifurcation point for \eqref{ODE} by 
Theorem \ref{thm} if we just can show that $\sfl(L,[a,b])\neq 0$.\\
In what follows, we restrict to the special case of positive definite matrices $A(x)$ in which our 
theory turns out to be particularly applicable. Let us introduce for $t\in[a,b]$ 
and $\lambda\in\mathbb{R}$ the following linear spaces
\begin{align*}
E(t_-,\lambda)&=\{u\in H_t:-(A(x)u'(x))'+S(x)u(x)=\lambda u(x),\, x\in[0,t]\}\\
E(t_+,\lambda)&=\{u\in H_t:-(A(x)u'(x))'+S(x)u(x)=\lambda u(x),\, x\in[t,1]\}
\end{align*}
as well as the non-negative integer
\[
\mu(t)=\sum_{\lambda<0}{\left(\dim E(t_-,\lambda)+\dim E(t_+,\lambda)+(\dim E(t_-,\lambda))\cdot(\dim E(t_+,\lambda))\right)}.
\]
Note that $\mu(t)< +\infty$ as $A(x)$ is positive definite for all $x\in I$.

\begin{prop}
Assume that the matrices $A(x)$, $x\in I$, are positive definite. If
\begin{enumerate}
	\item[(i)] $E(a_-,0)\cap E(a_+,0)=E(b_-,0)\cap E(b_+,0)=\{0\}$,
	\item[(ii)] $\mu(a)\neq\mu(b)$, 
\end{enumerate}
then there is a bifurcation point for \eqref{ODE}.
\end{prop}

\begin{proof}
It follows in our setting by \eqref{star}, \eqref{T} and \eqref{L} that
\[
\mu_{Morse}(L_t)\ =\ \sup\dim\{V\subset H_t:\, d^2_0\mathcal{J}[u,u]<0,\,\, u\in V\setminus\{0\}\}
\quad \mbox{for any $t \in I$}
\]
and so in view of Theorem \ref{thm} we need to show that:
\begin{enumerate}
	\item[(1)] the restrictions of $d^2_0J$ to $H_a$ and $H_b$ are non-degenerate,
	\item[(2)] $\mu_{Morse}(L_a) \ne \mu_{Morse}(L_b)$.	
\end{enumerate}
We note at first that if there exists $u\in H_t$ such that $d^2_0J[u,v]=0$ for all 
$v\in H_t$, then $u\in E(t_-,0)\cap E(t_+,0)$, which proves (1) by assumption (i).\\
Now, in order to show (2), we choose $\alpha>0$ such that the matrix 
$\alpha I_n+S(x)$ is positive definite for all $x\in[0,1]$,
where $I_n$ is the identity matrix on $\mathbb{R}^n$. Then, we get a scalar product on $H_t$ by
\[
\langle u,v\rangle_{t,\alpha}=\int^1_0{\langle A(x)u'(x),v'(x)\rangle\,dx}+
\int^1_0{\langle(\alpha I_n+S(x))u(x),v(x)\rangle\,dx},\quad u,v\in H_t,
\]
and by the Riesz representation theorem there exists a bounded operator $M$ on $H_t$ such that
\begin{align}\label{Riesz}
d^2_0J[u,v]=\langle Mu,v\rangle_{t,\alpha},\quad u,v\in H_t.
\end{align}
Hence $\mu_{Morse}(L_t)$ is the number of negative eigenvalues of $M$ counted with multiplicities. Now $Mu=\gamma u$ for some $\gamma<0$ if and only if
\begin{align*}
\langle Mu,v\rangle_{t,\alpha}&=\gamma \langle u,v\rangle_{t,\alpha}\\
&=\gamma\int^1_0{\langle A(x)u'(x),v'(x)\rangle\,dx}+\gamma\int^1_0{\langle S(x)u(x),v(x)\rangle\,dx}
+\gamma\alpha\int^1_0{\langle u(x),v(x)\rangle\,dx}
\end{align*}
for all $v\in H_t$. By \eqref{Riesz}, this is equivalent to
\[
- (A(x)u'(x))'+S(x)u(x)\ =\ \frac{\gamma\alpha}{1-\gamma}u(x),\quad x\in[0,t)\cup(t,1],
\]
and consequently, we see that 
\[
\mu_{Morse}(L_t)\ =\ \sum_{\lambda<0}{\dim\{u\in H_t:\, - (A(x)u'(x))'+S(x)u(x)=\lambda u(x),\, x\in[0,t)\cup(t,1]\}}.
\] 
Finally, there is a canonical isomorphism  $H^1_0([0,t],\mathbb{R}^n)\oplus H^1_0([t,1],\mathbb{R}^n)\rightarrow H_t$ 
which shows that the right hand side of the previous equality is indeed $\mu(t)$. 
\end{proof}
\noindent
Finally, let us mention that a related bifurcation problem is studied in \cite{AleIchDomain,AleIchBall,domainshrinking} and \cite{AleSmaleIndef}, where the authors consider the Dirichlet problem for elliptic partial differential equations 
\[
\left\{\begin{array}{ll}
- \Delta u + g(x,u) = 0 &\hbox{in $\Omega$}\\
u= 0 &\hbox{on $\partial\Omega$} 
\end{array}
\right.
\]
on a smooth bounded domain $\Omega\subset\mathbb{R}^N$ which is assumed to be 
star--shaped with respect to $0\in\mathbb{R}^N$. 
Denoting
\[
\Omega_t\ :=\ \{t x\in\mathbb{R}^N:\, x\in\Omega\}\subset\Omega \quad \hbox{for all $t\in[a,1]$,}
\] 
for some $0<a<1$, they study bifurcation of functionals on $H^1_0(\Omega,\mathbb{R})$ along the subspaces 
$\{H^1_0(\Omega_t,\mathbb{R})\}_{t\in[a,1]}$. 
However, our Theorems \ref{thm} and \ref{thmII} cannot be applied in this setting as the 
spaces $H^1_0(\Omega_t,\mathbb{R})$ do not vary continuously with respect to the metric of $\mathcal{G}(H^1_0(\Omega,\mathbb{R}))$. 
Indeed, if $0<s<t<1$, then there is a function $u\in H^1_0(\Omega_t,\mathbb{R})$ such that $\|u\|=1$ 
and with support in $\Omega_t\setminus\Omega_s$ (here, $\|\cdot\|$ is the standard norm in $H^1_0(\Omega,\mathbb{R})$). 
Consequently, $\langle u,v\rangle_{H^1_0(\Omega,\mathbb{R})}=0$ for all $v\in H^1_0(\Omega_s,\mathbb{R})$ and so
\[
\|P_{H^1_0(\Omega_t,\mathbb{R})}-P_{H^1_0(\Omega_s,\mathbb{R})}\|
\ \geq\ \|P_{H^1_0(\Omega_t,\mathbb{R})}u-P_{H^1_0(\Omega_s,\mathbb{R})}u\|\
=\ \|P_{H^1_0(\Omega_t,\mathbb{R})}u\|\ =\ \|u\| \ =\ 1,
\]    
which clearly contradicts the continuity.

\section*{Acknowledgements}
This work was commenced when the first author stayed  
at the {\sl Institut f\"ur Mathematik}, {\sl Humboldt Universit\"at zu Berlin},
and it was continued during a visit of the second author at the {\sl Dipartimento di Matematica}, {\sl Universit\`a degli Studi di Bari Aldo Moro}. The authors would like to thank both departments 
for their kind hospitality.


\vspace{1cm}
Anna Maria Candela\\
Dipartimento di Matematica\\
Universit\`a degli Studi di Bari ``Aldo Moro''\\
Campus Universitario\\
Via E. Orabona, 4\\
70125 Bari\\
Italy\\
E-mail: annamaria.candela@uniba.it

\newpage

Nils Waterstraat\\
School of Mathematics,\\
Statistics \& Actuarial Science\\
University of Kent\\
Canterbury\\
Kent CT2 7NF\\
UNITED KINGDOM\\
E-mail: n.waterstraat@kent.ac.uk
\end{document}